\documentclass[a4paper,10pt]{article}

\usepackage{amssymb}

\newtheorem{theorem}{Theorem}[section]
\newtheorem{lemma}[theorem]{Lemma}
\newtheorem{proposition}[theorem]{Proposition}

\newcommand{\proba}{{\mathbb P}}
\newcommand{\proof}{\noindent{\bf Proof}. }
\newcommand{\expectation}{{\mathbb E}}
\newcommand{\Oun}{\mathcal{O}(1)} 
\font\gfont=cmmi10 scaled \magstep2
\newcommand{\gchi}{\hbox{\gfont \char31}}
\newcommand{\ka}[1]{\kappa\big(\mathcal{E}_{#1},\mathcal{N}(0,\sigma^{2})\big)}

\newcommand{\keps}{[k^{{\scriptstyle \epsilon}/4}]}
\newcommand{\kumeps}{[k^{1-{\scriptstyle \epsilon}/4}]}

\begin{document}
\title{Almost Sure Central Limit Theorems\\
and the Erd{\"o}s-R{\'e}nyi law\\
for Expanding Maps of the Interval}
\bigskip

\author{J.-R. Chazottes\ \ \ and\ \ \ P. Collet\\
Centre de Physique Th{\'e}orique\\
CNRS UMR 7644\\
Ecole Polytechnique\\
F-91128 Palaiseau Cedex (France)\\
e-mail: jeanrene@cpht.polytechnique.fr\\
 collet@cpht.polytechnique.fr}

\maketitle

\begin{abstract}
For a large class of expanding maps of the interval, we prove 
that partial sums of Lipschitz observables satisfy an almost
sure central limit theorem (ASCLT). In fact, we provide a rate of convergence
in the Kantorovich distance. Maxima of partial sums are also
shown to obey an ASCLT. The
key-tool is an exponential inequality recently obtained.
Then we establish (optimal) almost-sure convergence rates for the supremum
of moving averages of Lipschitz observables (Erd{\"o}s-R{\'e}nyi
type law). This is done by refining the usual large deviations
estimates available for expanding maps of the interval.
We end up with an application to entropy estimation ASCLT's that refine
Shannon-McMillan-Breiman and Ornstein-Weiss theorems.
\end{abstract}

\section{Introduction and results} 

Almost sure central limit theorems were first derived for independent
sequences of random variables by Brosamler \cite{br}, Schatte
\cite{sc} and Fisher \cite{fi} (see also \cite{LP}). In its simplest form the result states that if
$X_{1},\;X_{2},\ldots$ is an i.i.d. sequence of real random variables 
with zero mean and unit variance (and satisfying some adequate moment
condition), if $(S_{n})$ is the sequence of
partial sums
$$
S_{n}=\sum_{j=1}^{n}X_{j}\;,
$$
then almost surely, for any real $x$
\begin{equation}\label{clt0}
\lim_{n\to\infty}\frac{1}{D_{n}}\sum_{k=1}^{n}\frac{1}{k}
\;\gchi_{(-\infty,x]}\left(\frac{S_{k}}{\sqrt{k}}\right)=\frac{1}{\sqrt{2\pi}}
\int_{-\infty}^{x}e^{-\frac{\xi^{2}}{2}}d\xi\;,
\end{equation}
where $D_{n}=\sum_{k=1}^{n}k^{-1}=\log n+\Oun$. 
Notice that this refines the following result by Erd{\"o}s-Hunt
\cite{EH} which can be seen as a discrete version of the strong arc
sine law: 
$$
\lim_{n\to\infty}\frac1{\log n} \sum_{k=1}^{n} \frac1{k} \gchi_{\{S_k>0\}}=\frac1{2}\quad\;\textup{a.s.}
$$
provided that the distribution of $X_1$ is symmetric.
The logarithmic average in (\ref{clt0}) may look surprising at first
glance. However, the Erd{\"o}s-Hunt result already implies that the
Ces{\'a}ro mean does not converge with probability one. We refer to
\cite{be,LP} for a more general discussion showing that the
logarithmic average is essentially the only one that works.

After the initial discovery, a large literature was published extending
the result in various directions. We refer to \cite{be}, \cite{bc} and
\cite{AW} for references and to \cite{lesigne} for a generic result
working for aperiodic dynamical systems.

Results for sequences of dependent random variables have also been
obtained under some mixing conditions (see \cite{be} and references
therein).  However dynamical
systems do not often satisfy these mixing conditions  with respect to
the most natural observables in their phase space. We will consider
below the family of dynamical systems given by piecewise expanding maps of the
interval.
We assume that the map $f$ is a piecewise monotonic transformation with
$b$ branches and we denote by
$(a_i)_{i=0}^{b}$ the corresponding subdivision of $[0,1]$.
We also assume that $f$ is monotonic and extends to a
$C^2$ map on each $\overline{A_i}=[a_i,a_{i+1}]$.
Finally we assume that the map is topologically mixing
and there is a constant $\eta>1$ and an integer $m$ such that
for any $x\in[0,1)$, $|(f^{m})'(x)|>\eta$.
Under these conditions this dynamical system 
has a unique absolutely continuous invariant measure $d\mu=h\;dx$ \cite{ly}
satisfying exponential decay of correlations for functions of bounded
variation \cite{hk}. For such observables one also has a central limit
theorem \cite{hk}. It is  therefore natural to investigate the truth of
almost sure central limit theorems in this class of systems which are
rather well understood.
We expect that the techniques developed below will
prove useful in more general situations.  
$\proba$ and $\expectation$ will refer to
the probability and the expectation with respect to $\mu$.

We will prove below a slightly stronger form of the almost sure central
limit theorem (\ref{clt0}) which is formulated using convergence in the
Kantorovich distance $\kappa$ (we refer to \cite{ra} for equivalent
definitions and properties of this distance).
 We will denote by $\mathcal{L}$ the set of
Lipschitz functions with Lipschitz constant equal to one and vanishing
at the origin. We recall that if $\nu_{1}$ and $\nu_{2}$ are two
probability measures, $\kappa$ is given by
\begin{equation}\label{dk}
\kappa\big(\nu_{1},\nu_{2}\big)=
\sup_{g\in\mathcal{L}}\int g(x)\;d\big(\nu_{1}-\nu_{2}\big)(x)\;.
\end{equation}
Note that since $\nu_{1}$ and $\nu_{2}$ are probability measures, in the
above integral we can replace $g$ by $g-g(0)$, or in other words there
is no restriction in assuming $g(0)=0$.
It is convenient  to define the sequence of  weighted empirical (random)
measures of the average by 
\begin{equation}\label{empisomme}
\mathcal{E}_{n}=\frac{1}{D_n}\sum_{k=1}^{n}\frac{1}{k}
\;\delta_{S_{k}/\sqrt k}\;,
\end{equation}
where $\delta_a$ is the Dirac measure at point $a$ of the real line.
We investigate the convergence of this random measures to a Gaussian distribution in the
Kantorovich metric.

We now state our first  main result.

\begin{theorem}\label{t1}
Consider  a piecewise $C^{2}$ expanding map $f$ of the interval which is
topologically mixing and let $d\mu=h\;dx$ be its unique absolutely
continuous invariant probability.  Consider  the sequence of random variables 
$X_{j}=u\circ f^{j}$ where $u$ is a Lipschitz function with zero $\mu$
average,   and assume that the  quantity $\sigma^{2}$ given by
\begin{equation}\label{variance}
\sigma^{2}=\int_{0}^{1} u^{2}d\mu+2\sum_{j=1}^{\infty}
\int_{0}^{1}u\circ f^{j}\;u\;d\mu
\end{equation}
is non zero. Then  Lebesgue almost surely
\begin{equation}\label{clt1}
\lim_{n\to\infty}\kappa\left(\mathcal{E}_{n},
\mathcal{N}\big(0,\sigma^{2}\big)\right)=0\;,
\end{equation}
where $\mathcal{N}\big(0,\sigma^{2}\big)$ is the Gaussian measure with
mean zero and variance $\sigma^{2}$.
More precisely, there is a constant $C_0>0$ such that, for all $n>3$,
$$
\proba\left(\sup_{j>n}\left( \frac{(\log j)^{1/3}\kappa\big({\cal E}_{j},{\cal
N}(0,\sigma^2)\big)}{\sqrt{\log\log j}}\right)>1\right)\le
$$
$$
\Oun\ e^{-C_0(\log n)^{1/3}}\;. 
$$
\end{theorem}

Notice that the last estimate in the theorem provides an upper bound
to the velocity of the approximation by the Gaussian measure.

It follows from the proof that
there exist constants $C_{1}>0$ and $C_{2}>0$ such that 
 for any $n>3$ and for any $t>0$
$$
\proba\left(\kappa\big({\cal E}_{n},{\cal
N}(0,\sigma^{2})\big)>t+C_{2}\frac{(\log\log n)^{1/2}}{(\log n)^{1/3}}\right)
\le 2\ e^{-C_{1}t^{2}D_{n}}\;.
$$

This is a kind of large deviation bound.  
In the independent case, and for the almost sure weak convergence
large deviation, estimates have been obtained in \cite{yor}.

It follows from the decay of correlations that the quantity
$\sigma^{2}$ is finite and non negative (see \cite{hk}).
We also recall that since $\mathcal{N}\big(0,\sigma^{2}\big)$ has a
bounded density, convergence in the Kantorovich metric implies
convergence in the Kolmogorov metric and in other topologies, see
e.g. \cite{gibbs} for a review.
We recall that the convergence in the Kantorovich metric follows from
the weak convergence and the convergence of the integral of the
function $g(x)=|x|$. However we will handle directly 
the Kantorovich metric in order to get estimates on the speed of convergence.

Define now the measure ${\cal G}(\sigma)$ by
$$
d{\cal G}(\sigma)(x):=\frac{\sqrt{2}}{\sqrt{\pi} \sigma}\vartheta(x)\
e^{-x^2/2\sigma}\ dx\,,
$$
where $\vartheta$ is the Heaviside function.
In the following theorem, $S^*_{n}=\sup_{k\leq n}S_{k}$.

\begin{theorem}\label{t2}
 Under the assumptions of theorem \ref{t1}, the sequence of random
 measures
\begin{equation}\label{max}
\mathcal{M}_{n}=\frac{1}{D_n}\sum_{k=1}^{n}\frac1{k}
\;\delta_{\frac{S^*_{k}}{\sqrt k}}\;,
\end{equation}
converges Lebesgue almost surely to ${\cal G}(\sigma)$
 in the Kantorovich metric.
\end{theorem}

The same result holds for the minima of partial sums $\inf_{k\leq n}S_{k}$.
Velocity of convergence estimates could also be obtained.

\bigskip

In order to state the  Erd{\"o}s-R{\'e}nyi Theorem, 
we recall the large deviation result for expanding maps of the interval
\cite{collet-lectures}.  Let $u$ be a function of bounded variation on $[0,1]$. 
Without loss of generality, we assume that
$\expectation(u)=0$. Suppose that $\sigma>0$ where $\sigma$ is defined
by formula (\ref{variance}). 

There exists $\beta_0>0$ such that for any
$\vert \beta\vert\leq \beta_0$ the following limit exists
$$
F(\beta)=\lim_{n\to\infty}\frac{1}{n}\log Z_{n}(\beta)
$$
where
$$
Z_{n}(\beta)=\int
e^{\beta\sum_{j=0}^{n-1}u(f^{j}(x))} d\mu(x)\;.
$$
For $\vert \beta\vert\leq \beta_0$ the function $F$ is analytic (see
\cite{broise}), and we denote by $\varphi$ 
its Legendre transform. Recall that the function $\varphi$ is analytic
for $\vert \alpha\vert\leq \alpha_0$, where $\alpha_0>0$ is small enough.
Moreover $\varphi$ is non negative
and strictly convex and its minimum, which is equal to zero, is attained at the $\mu$-expectation of
$u$, that is $0$.
In particular, for any fixed  $0<\alpha<\alpha_0$, and any
$\epsilon>0$, there is an integer $n_{0}=n_{0}(\epsilon)$ such that for
any integer $k>n_{0}$ we have
\begin{equation}\label{ep}
e^{-k(\varphi(\alpha)+\epsilon)}
\leq
\proba\left(\frac1{k}\sum_{j=0}^{k-1}u\circ f^{j}>\alpha\right)
\leq 
e^{-k(\varphi(\alpha)-\epsilon)}\;.
\end{equation}

We will need a much sharper result given in Appendix A. Let
$$
M_{k}=\sup_{0\le j\le \big[\exp(k\varphi(\alpha))\big]-k}
S_{k}\circ f^{j}\;.
$$

\begin{theorem}\label{t3}
Assume that $u$ is function of bounded variation of $[0,1]$. 
Then, under the hypotheses of Theorem \ref{t1},  
there is a number $\alpha^*>0$ {\rm (}$\alpha^*\leq \alpha_0${\rm )} such that for any $|\alpha|\leq
\alpha^*$, we have Lebesgue almost surely
$$
\limsup_{k\to\infty}\;\frac{M_{k}-k\alpha}{\log k}\leq \frac{1}{2\beta}
$$
and
$$
\liminf_{k\to\infty}\;\frac{M_{k}-k\alpha}{\log k}\geq -\frac{1}{2\beta}
$$
where $\beta=\varphi'(\alpha)$.
\end{theorem}

Note that this implies 
$$
\lim_{k\to\infty} \frac{M_k}{k}=\alpha\quad\textup{Lebesgue almost surely}\,.
$$

\noindent{\bf Remarks}. 

\noindent 1. We refer to \cite{deheuvels} for an optimal estimate in the case of independent
random variables. We notice that we get the same rate of convergence.
It follows from the estimate below that one can
derive estimates on the rate of convergence, however these estimates
depend on quantities like the derivative of the pressure function $F$ which are not easily
controlled in terms of data on the map $f$.  

\noindent 2. As observed by Comets \cite{comets} in a different context (Ising
model of Statistical Mechanics on a lattice), the 
Erd{\"o}s-R{\'e}nyi Theorem can be used as a statistical tool to determine
the large deviation function. As a first application, we recall that for
full Markov expanding maps of the interval, the
essential spectral radius of the transfer operator is related to the
function $F$ which is the inverse Legendre transform of $\varphi$ (see
\cite{colletisola}  and for extensions \cite{gundlachlatushkin}
and \cite{colleteckmann}). In other words, this Theorem provides a
statistical tool to estimate a generic lower bound on the decay of
correlations. Other applications of the Erd{\"o}s-R{\'e}nyi Theorem will
be discussed in forthcoming publications.

\noindent 3. It follows from results in \cite{buzzi} and \cite{liverani} 
that, since we assume that our map is topologically mixing, the density $h$ of the absolutely
continuous invariant measure $\mu$ is bounded below away from
$0$. Another sufficient condition is provided in \cite{barbour}.

\bigskip

The rest of the paper is organized as follows. In section 2 we prove
theorems \ref{t1} and \ref{t2}. We first derive an estimate on the
expectation of the Kantorovich distance. We next apply
an exponential  inequality proven in \cite{cms} which allows to
control the deviations of the Kantorovich distance from its 
expectation. The result follows from a Borel-Cantelli type
argument.
The Erd{\"o}s-R{\'e}nyi Theorem is proven in section 3 using a
precise large deviations estimate of independent interest which is
discussed in the appendix. In Section 4 we deduce ASCLT's refining
both Shannon-McMillan-Breiman and Ornstein-Weiss theorems. The latter
estimates the entropy of the system by the recurrence rate of typical
trajectories.

In the sequel, $C,C_1$, etc, will denote various constants whose value
may vary with the context. 

\section{Almost sure convergence for Lipschitz observables}

We will say that a real-valued function $K$ on $[0,1]^{n}$ is separately
Lipschitz in all its variables, if for any $1\le i\le n$ we have 
$$
L_{i}(K)=\sup_{x_{1}\cdots x_{n},\;y\in[0,1]}
$$
\begin{equation}\label{lip}
\frac{\bigg|K\big(x_{1},\cdots,x_{i-1},x_{i},x_{i+1},\cdots,x_{n}\big)
-K\big(x_{1},\cdots,x_{i-1},y,x_{i+1},\cdots,x_{n}\big)\bigg|}
{\big|x_{i}-y\big|} <\infty
\end{equation}
where in the above notation the indices out of range are absent. 

We recall some  inequalities  proven in \cite{cms}.

\begin{theorem}
There is a constant $C>0$
such that for any integer $n$ and for any 
real valued function   $K$ on  $[0,1]^{n}$ separately
Lipschitz in all its variables,  we have
\begin{equation}\label{ine}
\int e^{K(x,\cdots,f^{n-1}(x))-\expectation(K)}d\mu(x)
\le e^{C \sum_{j=1}^{n}L_j (K)^2}\;,
\end{equation}
where $\expectation(K)$ is the average
$$
\expectation(K)=\int K\big(x,\cdots,f^{n-1}(x)\big)d\mu(x)\;.
$$
We have also 
\begin{equation}\label{expe}
\proba\bigg(\big|K\big(x,\cdots,f^{n-1}(x)\big)-\expectation(K)\big|>t\bigg)
\le 2e^{-t^{2}/\big(4C\sum_{j=1}^{n}L_j(K)^2\big)}\;.
\end{equation}
\end{theorem}

Recall (see for example \cite{ra}) that the Kantorovich distance is
also given by
$$
\kappa\big({\cal E}_{n},{\cal N}(0,\sigma^{2})\big)=
\int_{-\infty}^{\infty}\big|{\cal F}_{n}(x)-F_{\sigma}(x)\big|\;dx
$$
where ${\cal F}_{n}$ and $F_{\sigma}$ are the laws of ${\cal E}_{n}$ and 
${\cal N}(0,\sigma^{2})$  respectively, or in other words
$$
{\cal F}_{n}(x)=\int_{-\infty}^{x} d{\cal E}_{n}=
\frac{1}{D_{n}}\sum_{1}^{n}\frac{1}{k}\vartheta\left(x-\frac{S_{k}}{\sqrt
k}\right)\;.
$$
($\vartheta$ is the Heaviside function.)

We will first prove that
$$
\lim_{n\to\infty}\expectation\left(\kappa\big({\cal E}_{n},{\cal
N}(0,\sigma^{2})\big)\right)=0\;. 
$$ 
In fact, we will prove a stronger result estimating the speed of convergence to
zero which is useful for the second part of Theorem \ref{t1}.

\begin{proposition}\label{propmoyenne}
Under the hypotheses of Theorem \ref{t1}, there is a positive constant
$C$ such that for any $n>3$ we have 
$$
\expectation\left(\kappa\big({\cal E}_{n},{\cal
N}(0,\sigma^{2})\big)\right)\le \frac{C(\log\log n)^{1/2}}{(\log n)^{1/3}}\;.
$$
\end{proposition}

\proof
First introduce the following notation for convenience:
$A_n=\hat{C}\sqrt{\log\log n}$,
where $\hat{C}$ is a large positive constant to be fixed later.
Introduce also the following sequence of functions:
$\psi_n(x)=\psi(x/A_n)$ where $\psi$ is a non-negative $C^2$
function bounded by one, equal to $0$ for $|x|>2$ and to $1$ for $|x|<1$.

It is easy to check that 
$$
\expectation\left(\kappa\big({\cal E}_{n},{\cal
N}(0,\sigma^{2})\big)\right)
\leq
(1+\Oun/A_n)\ \expectation\left(\sup_{g\in{\cal L}} \int_{A_n}^{\infty} g\ d({\cal E}_n-{\cal N}(0,\sigma^2))\right)
+
$$
$$
(1+\Oun/A_n)\ \expectation\left(\sup_{g\in{\cal L}} \int_{-\infty}^{-A_n} g\ d({\cal E}_n-{\cal N}(0,\sigma^2))\right)
+
$$
$$
\expectation\left(\sup_{g\in{\cal L}} \int g\psi_n d({\cal E}_n-{\cal N}(0,\sigma^2))\right)\;.
$$

Since $g(0)=0$ we have
$$
\expectation\left(\sup_{g\in{\cal L}} \int_{A_n}^{\infty} g\ d({\cal
    E}_n-{\cal N}(0,\sigma^2))\right)
\leq
\expectation\left( \int_{A_n}^{\infty} x\ d{\cal
    E}_n(x)\right) + \int_{A_n}^{\infty} x\ d{\cal N}(0,\sigma^2)(x)\;.
$$
The last integral is bounded above by $e^{-\Oun A_n^2}$ for $n$ large enough.
On the other hand, an integration by parts leads to
$$
\expectation\left(\int_{A_n}^{\infty} x\ d{\cal E}_n(x)\right)\leq
$$
$$
\frac{1}{D_{n}}\sum_{k=1}^{n}\frac{1}{k}
\int_{A_n}^{\infty}\proba(S_k>x\sqrt{k})\ dx +
\frac{1}{D_{n}}\sum_{k=1}^{n}\frac{1}{k}
A_n\proba(S_k>A_n\sqrt{k})\;. 
$$
Since the function
$$
K\big(x_{1},\cdots,x_{q}\big)=\left|\sum_{k=1}^{q}u(x_{k})\right|
$$ 
is separately Lipschitz with Lipschitz constants all equal to $L(u)$, it follows from
(\ref{expe}) that
there are two positive constants $C_{1}$ and $C_{2}$ such that for any integer $n$ and
any $t>0$
\begin{equation}\label{GDV}
\proba\bigg(\big|S_{q}\big|> t+\expectation\big(
\big|S_{q}\big|\big)\bigg)\le C_{1}e^{-C_{2}t^{2}/q}\;.
\end{equation}
It is proved in \cite{hk} that for any $q\geq 1$
\begin{equation}\label{GDVbis}
\expectation\big(|S_{q}|\big)\le
\expectation\big(S_{q}^{2}\big)^{1/2}\le\Oun \sqrt{q}\,.
\end{equation}
We now choose $\hat{C}$ large enough such
that for any $q\geq 1$ and $n\geq 3$ 
$$
A_n \sqrt{q} - \expectation(|S_q|) \geq \frac{A_n \sqrt{q}}{2}\,.
$$
Using (\ref{GDV}) with $t=x\sqrt{q} -\expectation(|S_q|)$
for $x\geq A_n$, and inequality (\ref{GDVbis}), we obtain
$$
A_n\proba(S_q>A_n\sqrt{q}) +
\int_{A_n}^{\infty}\proba(S_q>x\sqrt{q})\ dx\leq \Oun e^{-\Oun A_n^2}
$$
and therefore
$$
\expectation\left(\sup_{g\in{\cal L}} \int_{A_n}^{\infty} g\ d({\cal E}_n-{\cal N}(0,\sigma^2))\right)
\leq \Oun e^{-\Oun A_n^2}\;.
$$
Similarly
$$
\expectation\left(\sup_{g\in{\cal L}} \int_{-\infty}^{-A_n} g\ d({\cal E}_n-{\cal N}(0,\sigma^2))\right)
\leq \Oun e^{-\Oun A_n^2}\;.
$$

We now handle the middle range integral
$$
\expectation\left(\sup_{g\in{\cal L}} \int g\psi_n d({\cal E}_n-{\cal N}(0,\sigma^2))\right)\;.
$$

We introduce a sequence of finite sets ${\cal L}_{n}$ of Lipschitz
functions defined for $n>3$ and with Lipschitz constant at most $2$.
 This is a set of functions $g$ defined on
$[-2A_{n}-1,2A_{n}+1]$  
which are piecewise affine on each
interval between consecutive points of 
$({\mathbb Z}(\log n)^{-1/3})\cap[-2A_{n}-1,2A_{n}+1]$, 
which are zero at the origin and satisfy for any integer
$k\in (\log n)^{1/3}[-2A_{n}-1,2A_{n}+1]$  
$$
g\big((k+1)(\log n)^{-1/3}\big)-g\big(k(\log n)^{-1/3}\big)\in 
(\log n)^{-1/3}({\mathbb Z}\cap [-2,2])\;.
$$
Finally, one  takes $g$ constant on the
intervals
$$
\big[-2A_{n}-1,-[(2A_{n}+1)(\log n)^{1/3}](\log n)^{-1/3}\big]
$$
and
$$
\big[[(2A_{n}+1)(\log n)^{1/3}](\log n)^{-1/3}],2A_{n}+1\big]\;.
$$
It is easy to verify that for any $\ell\in{\cal L}$, there is a
$g\in{\cal L}_{n}$ such that 
$$
\sup_{x\in[-2A_{n}-1,2A_{n}+1]}|g(x)-\ell(x)|\le 3(\log n)^{-1/3}\;.
$$

\noindent By a simple computation one gets
\begin{equation}\label{cardin}
{\rm Card}\big({\cal L}_{n}\big)\le 5^{2(4A_{n}+3)(\log n)^{1/3}}\;.
\end{equation}
This estimate is related to  the $\epsilon$-entropy of ${\cal L}$, see \cite{kt}, \cite{lo}. 

We therefore have
$$
\expectation\left(\sup_{g\in {\cal L}}\int_{-\infty}^{+\infty}
\psi_n g\big(d{\cal E}_{n}-d{\cal N}(0,\sigma^2)\big)\right)\leq 
$$
$$
\expectation\left(\sup_{g\in {\cal L}_{n}}\int_{-\infty}^{+\infty}
\psi_n g\big(d{\cal E}_{n}
-d{\cal N}(0,\sigma^2)\big)\right)+\Oun A_{n}(\log
n)^{-1/3}\;. 
$$

We now use Pisier's inequality (see \cite{rio}) to estimate the integral
on the right hand side. We get
$$
\expectation\left(e^{(\log n)^{2/3}\sup_{g\in {\cal L}_{n}}
\int_{-\infty}^{+\infty} \psi_n g d({\cal
  E}_{n}-\mathcal{N}(0,\sigma^2))} \right)\le
$$
\begin{equation}\label{pisiermoyen}
\sum_{g\in {\cal L}_{n}} \expectation\left( e^{(\log n)^{2/3}
\int_{-\infty}^{+\infty} \psi_n g d({\cal E}_{n}-\mathcal{N}(0,\sigma^2))}\right)\;.
\end{equation}
The factor $(\log n)^{2/3}$ in the exponent will be convenient to 
balance later on the different bounds and, in particular, the `entropy'
contribution (\ref{cardin}).
We have
$$
\int_{-\infty}^{+\infty} \psi_n g\ d{\cal E}_{n}=
K_{n}\big(x,f(x),\cdots,f^{n-1}(x)\big)
$$
where 
$$
K_{n}\big(x_{1},\cdots,x_{n}\big)=
\frac{1}{D_{n}}\sum_{j=1}^{n}\frac{1}{j}
(\psi_n g)\left(\frac{\sum_{p=1}^{j}x_{p}}{\sqrt  j}\right)
$$
and this function is separately Lipschitz with
$$
L_{q}\big(K_{n}\big)\le \frac{\Oun}{D_{n}q^{1/2}}
$$
uniformly in $n$ and in $g\in{\cal L}_{n}$. We now apply Jensen's inequality to
the left hand side of (\ref{pisiermoyen}) and the exponential inequality
(\ref{ine}) to each term of the right hand side together with (\ref{cardin}) and we
obtain
$$
e^{(\log n)^{2/3}\expectation\left(\sup_{g\in {\cal L}_{n}}
\int_{-\infty}^{\infty} \psi_n g d({\cal E}_{n}-\mathcal{N}(0,\sigma^2))\right)}
$$
\begin{equation}\label{BE}
\le 5^{2(4A_{n}+3)(\log
n)^{1/3}} e^{\Oun (\log n)^{4/3} D_{n}^{-1}}
e^{(\log n)^{2/3}\sup_{g\in {\cal L}_{n}}\expectation\left(
\int_{-\infty}^{\infty} \psi_n g d({\cal E}_{n}-\mathcal{N}(0,\sigma^2))\right)}\;.
\end{equation}

We now take the logarithm on both sides and divide by $(\log
n)^{2/3}$. We obtain
$$
\expectation\left(\sup_{g\in {\cal L}_{n}}
\int_{-\infty}^{\infty} \psi_n g d({\cal
  E}_{n}-\mathcal{N}(0,\sigma^2))\right)\leq 
$$
$$
\frac{\Oun A_n}{(\log n)^{1/3}} +
\sup_{g\in {\cal L}_{n}}\expectation\left(
\int_{-\infty}^{\infty} \psi_n g d({\cal E}_{n}-\mathcal{N}(0,\sigma^2))\right)\,.
$$
We now estimate the last term in the right-hand side of this inequality.

Let $\nu_{k}$ be the law of $S_{k}/\sqrt k$. Since
$g\in{\mathcal L}_n$ is differentiable except in ${\mathbb Z}(\log n)^{-1/3}$
we have
$$
\int_{-\infty}^{\infty} \psi_n g d\nu_{k}\leq
\int_{-\infty}^{\infty}dx\, (\psi_n g)'(x) \nu_{k}\big((-\infty,x]\big)+
$$
$$
\sum_{|q|\leq (2A_{n}+1)(\log n)^{1/3}}\!\!\!\!\!\!\!\!
\Delta\nu_k (q(\log n)^{-1/3})\,
(\psi_n g)(q(\log n)^{-1/3}\big)\;,
$$
where
$$
\Delta\nu_k (s)=
\nu_{k}\big((-\infty,s]\big)-\nu_{k}\big((-\infty,s^-]\big)\;.
$$
To estimate $\Delta\nu_k$, we will use the Berry-Esseen inequality
proved in our context by A. Broise \cite{broise}). She indeed proved that
$$
\sup_{x\in\mathbb{R}}\left|\expectation\left(\vartheta\left(x-\frac{S_{k}}{\sqrt{k}}\right)
\right)-F_{\sigma}(x)\right|\le\Oun k^{-1/2}
$$
uniformly in $k\ge1$.
Since the Gaussian has a
bounded density, we have uniformly in $k$
$$
\sup_{|q|\leq (2A_{n}+1)(\log n)^{1/3}}\!\!\!\!\!
\Delta\nu_k (q(\log n)^{-1/3})
\le
\Oun k^{-1/2}
$$
therefore uniformly in $n$ and $k$
$$
\int_{-\infty}^{\infty} \psi_n g\ d\nu_{k}=\int_{-\infty}^{\infty}
\psi_n g\ d{\cal
N}(0,\sigma^2) +\Oun A_{n}(\log n)^{1/3} k^{-1/2}\;.
$$
From (\ref{empisomme}) we have
$$
\expectation\left(
\int_{-\infty}^{\infty} \psi_n g d{\cal E}_{n}\right)= 
\frac{1}{D_n} \sum_{k=1}^{n} \frac{1}{k} 
\int_{-\infty}^{\infty} \psi_n g\ d\nu_{k}\,.
$$
Then it follows easily that uniformly in $n\ge 1$ and $g\in{\cal L}_{n}$ we have
$$
\expectation\left(
\int_{-\infty}^{\infty} \psi_n g d{\cal E}_{n}\right)= 
\int_{-\infty}^{\infty} \psi_n g d{\cal N}(0,\sigma^2)+\Oun A_{n}(\log n)^{-2/3}\;. 
$$
Using this estimate in (\ref{BE}), the proposition follows.

\bigskip

We now handle the proof of Theorem \ref{t1}.
In order to apply the exponential inequality, we consider the sequence
of functions
$$
K_{n}\big(x_{1},\ldots,x_{n}\big)=
\int_{-\infty}^{\infty}dx\left|\frac{1}{D_{n}}\sum_{k=1}^{n}\frac{1}{k}
\vartheta\left(x-\frac{\sum_{l=1}^{k}x_{l}}{\sqrt k}\right)
-F_{\sigma}(x)\right|\;.
$$
After an easy computation one gets for any $1\le p\le n$
$$
L_{p}\big(K_{n}\big)\le \frac{{\cal O}(1)}{D_{n}p^{1/2}}\;.
$$
We can now apply the exponential inequality (\ref{ine}) to get the
existence of a constant $C_{1}>0$ such that for any $n$ and for any $t>0$
$$
\proba\left(\bigg|\kappa\big({\cal E}_{n},{\cal N}(0,\sigma^{2})\big)
-\expectation\left(\kappa\big({\cal E}_{n},{\cal
N}(0,\sigma^{2})\big)\right)\bigg|>t\right)\le 2e^{-C_{1}t^{2}D_{n}}\;.
$$

This implies the first part of Theorem \ref{t1}. 

We define the sequence $(n_k)$ by
$$
n_{k}=e^{k^{3}}\;.
$$
We conclude from the above estimate and Proposition \ref{propmoyenne}
that for a positive constant $\Theta$ large enough
$$
\sum_{k}
\proba\left(\ka{n_k}>\frac{\Theta(\log\log n_k)^{1/2}}{(\log n_k)^{1/3}}\right)
<\infty
$$
which implies by the Borel-Cantelli lemma that the sequence
$\ka{n_k}$ converges to zero almost
surely.
More precisely
$$
\sum_{k\geq (\log n)^{1/2}}
\proba\left(\ka{n_k}>\frac{\Theta(\log\log n_k)^{1/2}}{(\log n_k)^{1/3}}\right)
\leq
$$
$$
\Oun \sum_{\ell\geq (\log n)^{1/3}} e^{-\Oun \ell}\leq \Oun
e^{-\Oun (\log n)^{1/3}}\;. 
$$
From (\ref{empisomme}) it follows that for
$n>n_{k}$
$$
\mathcal{E}_n = \mathcal{E}_{n_k} -
\frac{D_{n}-D_{n_k}}{D_{n}}\ \mathcal{E}_{n_k} 
+ \frac{1}{D_n} \sum_{j=n_k +1}^{n}\frac{1}{j} \delta_{S_j /\sqrt{j}}\,. 
$$
This implies
$$
\left|\ka{n}-\ka{n_{k}}\right|\leq
$$
$$
\frac{D_{n}-D_{n_k}}{D_{n}}\ka{n_{k}}+
\sup_{g\in \mathcal L}\frac{1}{D_n} \sum_{j=n_{k}+1}^{n}\frac{1}{j}\left[
g\left(\frac{S_{j}}{\sqrt
j}\right)-\int g\ d\mathcal{N}(0,\sigma^{2})
\right].
$$
From now on, we assume that $n\le n_{k+1}$, then
the first term tends to zero almost surely by our previous estimates
and is more precisely $\mathcal{O}((\log n_k)^{-1/3})$. We
now prove that the second term tends to zero almost surely. 
We have
$$
\sup_{g\in \mathcal L}\frac{1}{D_n} \sum_{j=n_{k}+1}^{n}\frac{1}{j}\left[
g\left(\frac{S_{j}}{\sqrt
j}\right)-\int g\ d\mathcal{N}(0,\sigma^{2})
\right]
$$
$$
\le  \frac{1}{D_{n}}\sum_{j=n_{k}+1}^{n}\frac{1}{j}\left[
\frac{\big|S_{j}\big|}{\sqrt j}+\int |x|\ d\mathcal{N}(0,\sigma^{2})(x)
\right]
$$
$$
\le  \frac{1}{D_{n_{k}}}\sum_{j=n_{k}+1}^{n_{k+1}}\frac{1}{j}\left[
\frac{\big|S_{j}\big|}{\sqrt j}+\int |x|\ d\mathcal{N}(0,\sigma^{2})(x)
\right]\;.
$$
We have for any $k$
$$
\frac{1}{D_{n_{k}}}\sum_{j=n_{k}+1}^{n_{k+1}}\frac{1}{j}
\int |x| d\mathcal{N}(0,\sigma^{2})(x)\leq \Oun \frac{\log n_{k+1}  -\log
  n_k }{\log n_k }
$$
Using (\ref{GDVbis}), we get
$$
\expectation\left( \frac{1}{D_{n_k}}\sum_{j=n_{k}+1}^{n_{k+1}}\frac{1}{j}
\frac{\big|S_{j}\big|}{\sqrt j} \right)
\leq 
\Oun \frac{\log n_{k+1}  -\log n_k }{\log n_k }\;.
$$
We now observe that the function
$$
K(x_1,...,x_{n_{k+1}})=
\frac{1}{D_{n_{k}}}\sum_{j=n_{k}+1}^{n_{k+1}}\frac{1}{j}
\frac{\big|\sum_{\ell=1}^{j}u(x_\ell)\big|}{\sqrt j}
$$
is separately Lipschitz with the following estimates, for $k>4$:
$$
L_q(K)\leq\cases{\frac{\Oun}{\sqrt{n_k}D_{n_k}} & for
  $1\leq q<n_k + 1$ \cr
\frac{4}{\sqrt{q} D_{n_k}} & for $n_{k} + 1 \leq q\leq n_{k+1}$.\cr
}
$$
Therefore 
$$
\sum_{q=1}^{n_{k+1}} L_q(K)^2\leq  \Oun \frac{D_{n_{k+1}} -
    D_{n_k}}{D_{n_k}^2}\;.
$$
Using (\ref{expe}) and the choice of $n_k$, one easily gets for $k>1$
$$
\proba\left(\sup_{\ell\geq k}\left((\log n_\ell)^{1/3}
\frac{1}{D_{n_{\ell}}}\sum_{j=n_{\ell}+1}^{n_{\ell+1}}\frac{1}{j}
\frac{\big|S_j\big|}{\sqrt j}\right) >(\log\log n_\ell)^{1/2}\right)\le
$$
$$
\Oun \sum_{\ell\geq k} e^{-C_0 (\log n_\ell)^{1/3}}
$$
for some positive constant $C_0$.

Collecting all the above estimates finishes the proof of Theorem \ref{t1}.

\bigskip

The proof of Theorem \ref{t2} is similar to the proof of Theorem
\ref{t1}. Using the exponential inequality (\ref{expe}), it follows as before that
$$
\proba\left(\left|\kappa\big({\cal M}_{n},{\cal G}(s)\big)-
\expectation\big(\kappa\big({\cal M}_{n},{\cal
G}(s)\big)\big)\right|>t\right) \le 2 e^{-\Oun t^{2}/D_{n}}\;.
$$
We have also
$$
\proba\bigg(\big|S_{n}^{*}\big|>u+\expectation\big(\big|S_{n}^{*}\big|\big)
\bigg)\le \Oun e^{-\Oun u^{2}/n}\;.
$$
Let us now estimate the expectation. 
First note that 
$
\expectation\left(|S_n^*|\right)\leq \expectation\left(\sup_{1\leq
    j\leq n} |S_j|\right)
$.
We use Pisier's inequality \cite{rio} to get
$$
\expectation
\left(e^{\sup_{1\leq j\leq n} |S_j|\sqrt{\log n}/\sqrt{n}}\right)
\leq
\sum_{j=1}^n 
\expectation\left(e^{|S_j|\sqrt{\log n}/\sqrt{n}}\right)\;.
$$
Taking logarithm of both sides, using Jensen's inequality in the left hand side and inequality
(\ref{ine}) in the right hand side we get
$$
\expectation\big(\big|S_{n}^{*}\big|\big)
\le \Oun \sqrt{n\log n}\,.
$$
Therefore, for any $\eta>0$ we have
$$
\sum_{n=1}^{\infty}\proba\bigg(\big|S_{n}^{*}\big|>\sqrt n\, (\log
n)^{1/2+\eta}\bigg)<\infty\;. 
$$

We now deal with the expectation.

\begin{proposition}
Under the assumptions of Theorem \ref{t1}, we have
$$
\lim_{n\to\infty}\expectation(\kappa({\cal M}_n,{\cal G}(\sigma)))=0\;.
$$
\end{proposition}

\proof
As before we have for any $A>0$ large enough and uniformly in $n$ 
$$
\expectation\left(\int_{|x|>A}dx
\left|\frac{1}{D_{n}}\sum_{1}^{n}\frac{1}{k}\vartheta\left(x-\frac{S^*_{k}}{\sqrt k}\right)
-{\cal G}(\sigma)(x)\right|\right)
\leq e^{-\Oun A^2}\;.
$$

The middle integral is treated as before except for the term
$\expectation\left(\int \psi_A g\ d{\cal M}_n\right)$ appearing in the
analog of (\ref{BE}). 

It is enough to show that for any $g\in {\cal L}_n$
$$
\lim_{k\to\infty}\expectation\left((\psi_A g) (S^*_k/\sqrt{k})
\right)=\int \psi_A g d{\cal G}(\sigma)\;.
$$
At this point we recall the almost sure invariance principle (ASIP)
\cite{de}, \cite{hk}.
There is a positive number $\delta$ and an enriched probability space $\Omega$ carrying also a Brownian
motion $(B_t)$ such that 
$$
\tilde{S}_k - \sigma B_k = o(k^{-\delta+1/2})\quad\textup{eventually almost
  surely}
$$
where the sequence $(\tilde{S}_k)$ has the same joint distribution 
as $(S_k)$. In other words there is an almost-surely finite
integer-valued random variable $N(\omega)$, where $\omega\in\Omega$, such that, for any
$j>N(\omega)$
$$
| \tilde{S}_j (\omega) - \sigma B_j(\omega)| \leq j^{-\delta+1/2}\;.
$$
To derive the ASIP for $(\tilde{S}_k^*)=(\sup_{j\leq k} \tilde{S}_k)$, notice that there are two
cases for each $k$, namely the supremum is attained for an index
$j_k\leq N(\omega)$ or for an index $j_k>N(\omega)$. 
We claim that for almost all $\omega$ the first case can occur at
most for finitely many $k$'s. Indeed, if it was not the case, since
$(\tilde{S}_k^*)$ is non decreasing, it would imply that this sequence
is bounded. However, by the law of the iterated logarithm (for $B_k$)
the sequence $(\tilde{S}_k)$ diverges almost surely. In other words, 
the index for which the maximum is attained in
$\tilde{S}_k^*$ is eventually almost surely larger than $N(\omega)$.
The same argument holds for the sequence $(B_j)$.
We conclude that the ASIP also holds for the sequence
$(\tilde{S}^*_k)$. In particular
$$
\frac{\tilde{S}_k^*}{\sqrt{k}} -\frac{\sigma B_k^*}{\sqrt{k}}\to
0\quad\textup{almost surely}\;.
$$
Since $\sup_{x\in\mathbb{R}}|(g\psi_A)(x)|\leq 2 A$, we can use the dominated convergence theorem to conclude that for
any $A>0$   
$$
\lim_{k\to\infty}\left(\expectation\big((\psi_A
  g)(S_k^*/\sqrt{k})\big)-
\expectation\big((\psi_A g)(\sigma B_k^*/\sqrt{k})\big) \right)=0\;.
$$
Define $B_k^{**}=\sup_{1\leq j\leq k} B_{j/k}$. Now observe that by
rescaling we get
$$
\expectation\big((\psi_A g)(\sigma B_k^*/\sqrt{k})\big)=
\expectation\big((\psi_A g)(\sigma B_k^{**})\big)\;.
$$
Since trajectories of the Brownian motion are almost surely H{\"o}lder
continuous \cite{gk}, we obtain
$$
B_k^{**} - \sup_{0\leq t\leq 1} B_t \to 0 \quad\textup{almost surely}\;.
$$
Using again the dominated convergence theorem we get
$$
\lim_{k\to\infty}\expectation\big((\psi_A g)(\sigma B_k^{**})\big)=
\expectation\big((\psi_A g)(\sigma \sup_{0\leq t\leq 1} B_t )\big)
$$
The proof of Theorem \ref{t2} is complete by letting $A$ tend to
infinity and by using the explicit expression for the law
of $\sup_{0\leq t\leq 1} B_t$, see \cite{gk}.


\section{Proof of Theorem \protect{\ref{t3}}}

From Lemma \ref{estime}, stationarity  and the differentiability of
$\varphi$,  we  have
$$
\proba\big(M_{k}>k\alpha+x\log k\big)
\le\sum_{j=0}^{[\exp(k\varphi(\alpha)]-k}
\proba\big(S_{k}\circ f^{j}>k\alpha+x\log k\big)
$$
$$
\le \Oun e^{k\varphi(\alpha)}k^{-1/2} e^{-k\varphi(\alpha+(x\log k)/k))}\le \Oun 
k^{-1/2}e^{-x\varphi'(\alpha)\log k}   
$$
and this is summable over $k$ if $x>1/(2\varphi'(\alpha))=1/(2\beta)$.
We are done with the upper bound.

\bigskip

For an integer  $r$ to be chosen later on independently of
$k$ and $0<\epsilon<1$ also independent of $k$, we define a new quantity
$$
\tilde M_{k}=\sup_{0\le j\le \big[(\exp(k\varphi(\alpha))-k)/(rk)]}\;
\sup_{0\le l\le \kumeps}
S_{k}\circ f^{jrk+l\keps}\;.
$$
The gaps of size $rk$ in the indices will allow us later on to use the
decay of correlations.  The small gap $\keps$ is chosen for
convenience. It can be reduced to a large enough constant times $\log k$. Its
role is to ensure that the probability of having simultaneously 
$S_{k}(x)>k\alpha$ and 
$S_{k}\circ f^{\keps}(x)>k\alpha$   is small enough. We have of
course

$$
\tilde M_{k}\le M_{k}
$$
and therefore for any $\lambda$
$$
\proba\big(M_{k}<\lambda)\le \proba\big(\tilde M_{k}<\lambda)\;.
$$
We will now estimate this last quantity with $\lambda=k\alpha
-(1+\epsilon)\log k/(2\beta)$.
Let 
$$
A=\big\{x\;~:\; S_{k}(x)<k\alpha -(1+\epsilon)(\log k)/(2\beta)\big\}\;.
$$
We have
$$
\proba\big(\tilde M_{k}<k\alpha -(1+\epsilon)(\log k)/(2\beta)
\big)=
$$
$$
\expectation\left(
\prod_{{\scriptstyle 0\le j\le \big[(\exp(k\varphi(\alpha))-k)/(rk)]}
\atop 0\le l\le \kumeps}\gchi_{A}\circ
f^{jrk+l \keps}\right)\;. 
$$
In order to estimate this quantity, we will construct an upper
bound as follows. Let $\psi$ be the non negative Lipschitz function defined by
$$
\psi(x)=\cases{1&for $x<0$\cr
1-x&for $0\le x\le 1$\cr
0&for $1\le x$.\cr}
$$
We have obviously
\begin{equation}\label{psisup}
\gchi_{A}(x)\le\psi\big(S_{k}(x)-k\alpha +(1+\epsilon)(\log
k)/(2\beta)\big)\le\gchi_{B}(x)
\end{equation}
where
$$
B=\big\{x\;:\; S_{k}(x)<k\alpha -(1+\epsilon)(\log k)/(2\beta)+1\big\}\;.
$$
Therefore 
$$
\proba\big(\tilde M_{k}<k\alpha -(1+\epsilon)(\log k)/(2\beta)
\big)
$$
$$
\le \expectation\left(
\prod_{{\scriptstyle 0\le j\le \big[(\exp(k\varphi(\alpha))-k)/(rk)]}
\atop 0\le l\le \kumeps}\psi\circ
f^{jrk+l \keps}\right)\;. 
$$

It is easy to verify  that if $v$ is a function of bounded variation and
$w$ is a Lipschitz function with Lipschitz constant $K_{w}$ we have with
$\vee$ denoting the total variation
$$
\vee (w\circ v)\le K_{w}\vee v\;.
$$
We recall (see \cite{hk}) that under our hypothesis on the
transformation $f$, there is a positive constant $A>2$ such that for any
integer $q$ we have
$$
\vee f^{q}\le A^{q}\;.
$$
Therefore, if we define the function $g_{1}$ by
$$
g_{1}(x)=\prod_{0\le l\le \kumeps}
\psi\left(S_{k}\bigg(f^{l\keps}(x)\bigg)-k\alpha +(1+\epsilon)(\log
k)/(2\beta)\right)\;,
$$
one gets easily
$$
\vee g_{1}
\le k^{2}  A^{2k}\;.
$$
We recall (see \cite{hk}) that there exist two
positive constants $C$ and $\rho<1$ such that if $g_{1}$ is a function
of bounded variation and $g_{2}$ is integrable with respect to the
Lebesgue measure (and hence also with respect to $\mu$), we
have for any integer $p$ 
$$ 
\left|\int g_{1}\;g_{2}\circ f^{p}d\mu-\int g_{1}d\mu\int g_{2}d\mu\right|
\le C\rho^{p}\Big(\vee g_{1}+\int | g_{1} | \ dx\Big)\ \int | g_{2} | \ dx\;.
$$
We now apply this inequality with $p=(r-1)k$ 
and 
$$
g_{2}(x)=g_{2,q}(x)=\prod_{\scriptstyle 0\le j\le q-1
\atop 0\le l\le \kumeps}
\psi\left(S_{k}\big(f^{jrk+l\keps }(x)\big)
-k\alpha +(1+\epsilon)(\log k)/(2\beta)\right)\;,
$$
$q$ being an integer. We obtain the recursive bound for any $q\ge 1$
$$
\int d\mu(x)\,\prod_{\scriptstyle 0\le j\le q
\atop 0\le l\le \kumeps}
\psi\left(S_{k}\big(f^{jrk+l\keps}(x)\big)
-k\alpha +(1+\epsilon)(\log k)/(2\beta)\right)
$$
$$
\le \left(\int g_{1}d\mu+Ck^2A^{2k}\rho^{(r-1)k}\right)
$$
$$
\times \int d\mu(x)\,\prod_{\scriptstyle 0\le j\le q-1
\atop 0\le l\le \kumeps}
\psi\left(S_{k}\big(f^{jrk+l\keps}(x)\big)
-k\alpha +(1+\epsilon)(\log k)/(2\beta)\right)\;.
$$
We now choose $r$ (depending on $\alpha$) 
such that $2\log A+(r-1)\log\rho<-\varphi(\alpha)$. This
implies for $k$ large enough
$$
\proba\big(M_{k}<k\alpha -(1+\epsilon)(\log k)/(2\beta)
\big)
$$
\begin{equation}\label{chat}
\le
\left(\int g_{1}d\mu +\Oun
e^{-k\varphi(\alpha)-k\delta}\right)^{\big[
(\exp(k\varphi(\alpha))-k)/(rk)\big] -1} 
\end{equation}
for some $\delta>0$ independent of $k$.

We now have to estimate the integral of $g_{1}$ from above. We have of
course from (\ref{psisup})
$$
\int g_{1}(x)d\mu \le \int d\mu \prod_{0\le l\le \kumeps}
\gchi_{B}\circ f^{l\keps}
$$
$$
=\int d\mu
\prod_{0\le l\le \kumeps}
\left(1-\gchi_{B^c}\circ f^{l\keps}\right)
\;.
$$
Using Bonferoni's inequality, we get
$$
\int g_{1}(x)d\mu \le 1-\sum_{0\le l\le \kumeps}\int
d\mu \, \gchi_{B^c}\circ f^{l\keps}
$$
\begin{equation}\label{bonferoni}
+\sum_{0\le m\neq n \le \kumeps}\int
d\mu \, \gchi_{B^c}\circ f^{n\keps} \gchi_{B^c}\circ f^{m\keps}\;.
\end{equation}

We observe that for fixed $0<\epsilon<1$ and for $k$ large enough, we
have from Lemma (\ref{estime})
$$
\mu(B^c)\ge \Oun\
k^{\frac{\epsilon}{2}}\
e^{-k\varphi(\alpha)} \;.
$$
Therefore using also the invariance of the measure $\mu$, the opposite
of the second  term in the right hand side is bounded below by
\begin{equation}\label{souris}
(\kumeps+1)\mu(B^c)\ge\Oun \ k^{1+\epsilon/4} e^{-k\varphi(\alpha)}\;.
\end{equation}
We now have to estimate the last term on the right hand side of
(\ref{bonferoni})  and show in
particular that it is much smaller than the modulus of the second term
(for large $k$). This estimate is provided in the appendix by
Lemma \ref{decouple}. 
Collecting all the bounds, namely (\ref{chat}), (\ref{bonferoni}) and (\ref{souris}), we finally get
$$
\proba\big(M_{k}<k\alpha -(1+\epsilon)(\log k)/(2\beta)
\big)\le \Oun e^{-C \keps}
$$
where $C>0$ is some positive constant.
The right hand side of this estimate 
 is summable in $k$ for any $\epsilon>0$ and the result
follows using the Borel-Cantelli Lemma.

\section{Applications to entropy estimation}

We define $n$-cylinder sets as usual:
$A_{i_1}^{i_n}:= A_{i_1}\cap f^{-1}A_{i_2}\cap\cdots\cap
f^{-n+1}A_{i_n}$,
where the $A_i$'s are the monotonicity/regularity intervals of the map
$f$. Denote by
$\mathcal{P}_n$ the set of $n$-cylinders. For all $x\in[0,1]$ which is not the $n$-th preimage
of a discontinuity point, there is a unique $n$-cylinder containing $x$, denoted by $\mathcal{P}_n(x)$.
Since $f$ is expanding, the partition of $[0,1]$ into the sets $]a_i,a_{i+1}]$ generates the Borel
$\sigma$-algebra.

We assume that $\log|f'|$ is a Lipschitz function. Throughout this section,
the observable $u$ will be $\log|f'|$.
Recall that by Rokhlin formula \cite{ledrap}, $h_\mu(f)=\int \log |f'|\ d\mu$.

The following theorem is a refinement of Shannon-McMillan-Breiman
Theorem (see also \cite{ibragimov},
\cite{ziemian} for the
usual central limit theorem). Let us define the sequence of  weighted empirical (random)
measures of minus the logarithm  of the $\mu$-measure of cylinders by 
$$
\mathcal{SMB}_{n}=\frac{1}{D_n}\sum_{k=1}^{n}\frac{1}{k}
\;\delta_{(-\log\mu(\mathcal{P}_k(\cdot))-kh)/\sqrt k}
$$
where $h=h_\mu(f)$.

\begin{theorem}\label{asclt-smb}
Under the hypotheses of Theorem \ref{t1} and assuming that the
function $\log|f'|$ is Lipschitz, 
\begin{equation}
\lim_{n\to\infty}\kappa\left(\mathcal{SMB}_{n},
\mathcal{N}\big(0,\sigma^{2}\big)\right)=0\;,
\end{equation}
Lebesgue almost surely.
\end{theorem}

\proof
By Theorem \ref{t1} applied to $u=\log|f'|-h$ and the triangle
inequality,
we only have to prove that $\kappa(\mathcal{SMB}_n,\mathcal{E}_n)\to 0$
as $n$ goes to infinity, Lebesgue almost surely, where $\mathcal{E}_n$
is given by (\ref{empisomme}).

We have the following strong approximation:
for any $0<\varrho<1$ and for any $x$ in a set of measure $\geq 1-\varrho$, there is
an integer $N=N(\varrho)$ such that for any $k>N$
\begin{equation}\label{Smu}
S_k(x) -\mathcal{O}(1) \leq -\log \mu(\mathcal{P}_k(x)) \leq S_k(x) +\mathcal{O}(1)\,.
\end{equation}
Indeed, let
$$
B(k,C)=\left\{P\in\mathcal{P}_k: \forall x\in P, \frac1{C}\leq \frac{\mu(P)}{\exp -S_k(x)}\leq C\right\}\,.
$$
Then, by Lemma 22 in \cite{paccaut-PhD}, given any $0<\varrho<1$, there exists $\Gamma_\varrho>0$ and $N_\varrho>0$ such 
that
$$
\mu\left(\bigcap_{n>N_{\varrho}}\bigcup_{P\in B(n,\Gamma_\varrho)}P\right)\geq 1-\varrho\,.
$$

As before, we can write this distance as follows:

$$
\kappa(\mathcal{SMB}_n,\mathcal{E}_n)\leq
$$
$$
\frac1{D_n}\sum_{k=1}^n \frac1{k}\int_{-\infty}^{+\infty} d\xi\
\Big\vert \vartheta\left(\xi -\frac{S_k-kh}{\sqrt{k}}\right)-
\vartheta\left(\xi -\frac{-\log\mu(\mathcal{P}_k(\cdot))-kh}{\sqrt{k}}\right)\Big\vert\, .
$$

It is now straightforward to get using (\ref{Smu})
$$
\kappa(\mathcal{SMB}_n,\mathcal{E}_n)\leq \frac1{D_n}\sum_{k=1}^n \frac{\mathcal{O}(1)}{k^{3/2}}
\leq \frac{\mathcal{O}(1)}{D_n}\to 0\quad\textup{as}\,n\to\infty
$$
on a set of measure $\geq 1-\varrho$.
Since $\varrho$ can be chosen arbitrarily small, the proof is finished.

\bigskip

We now define the return time of a point $x$ to its $n$-cylinder by
$$
R_n(x)=\inf\{k>0: f^k(x)\in \mathcal{P}_n(x)\}\,.
$$

The following result is due to Ornstein and Weiss \cite{ow}:
$$
\lim_{n\to\infty}\frac1{n}\log R_n(x)= h_\mu(f)\quad\textup{for}\,\mu\,\textup{almost all}\,x\,.
$$
This result needs only the ergodicity of $\mu$ to hold.

We define the sequence of  weighted empirical (random)
measures of the logarithm of the return times by 
$$
\mathcal{OW}_{n}=\frac{1}{D_n}\sum_{k=1}^{n}\frac{1}{k}
\;\delta_{(\log R_k-kh) /\sqrt k}\;.
$$

\begin{theorem}\label{asclt-ow}
Under the hypotheses of Theorem \ref{t1}, and
assuming that the function $\log|f'|$ is Lipschitz, we have
\begin{equation}\label{clt2}
\lim_{n\to\infty}\kappa\left(\mathcal{OW}_{n},
\mathcal{N}\big(0,\sigma^{2}\big)\right)=0\;,
\end{equation}
Lebesgue almost surely.
\end{theorem}

The main point is to get a sufficiently strong approximation of $\log R_k$ by
$S_k$ that holds eventually almost surely. This is the subject of Lemma \ref{lemma}
which is stated and proved in Appendix \ref{B}.

\proof
By the triangle inequality, it is enough to prove that $\kappa(\mathcal{OW}_n,\mathcal{E}_n)\to 0$
as $n$ goes to infinity, Lebesgue almost surely.

By Lemma \ref{lemma} and (\ref{Smu}), we obtain the following strong approximation:
for any $0<\varrho<1$ and for any $x$ in a set of measure $\geq 1-\varrho$, there is
an integer $N=N(x,\varrho)$ such that for all $n>N$
\begin{equation}\label{SR}
S_k(x) -\mathcal{O}(1)\log k \leq \log R_k(x) \leq S_k(x) +\mathcal{O}(1)\log k\,.
\end{equation}

We can write this distance as follows (as we have already done before):
$$
\kappa(\mathcal{OW}_n,\mathcal{E}_n)\leq 
\frac1{D_n}\sum_{k=1}^n \frac1{k}\int_{-\infty}^{+\infty} d\xi\
\Big\vert \vartheta\left(\xi -\frac{S_k-kh}{\sqrt{k}}\right)-
\vartheta\left(\xi -\frac{\log R_k -kh}{\sqrt{k}}\right)\Big\vert\, .
$$

It is straightforward to get using (\ref{SR})
$$
\kappa(\mathcal{OW}_n,\mathcal{E}_n)\leq \frac1{D_n}\sum_{k=1}^n \frac{\mathcal{O}(1)\log k}{k^{3/2}}
\leq \frac{\mathcal{O}(1)}{D_n}\to 0\quad\textup{as}\,n\to\infty
$$
on a set of measure $\geq 1-\varrho$.
Since $\varrho$ can be chosen arbitrarily small, the proof is finished.

We recall that lognormal fluctuations of the return time $R_n$ have been studied in \cite{cgs,konto,paccaut}.

\appendix
\section{Appendix: Sharp large deviation estimates.}

In this section we will establish a sharp estimate on the probability of
large deviations. Similar results have been obtained for independent
random variables and Markov chains. We refer to \cite{petrovbis}, 
\cite{INN}, \cite{Iltis} and \cite{ney} for a detailed discussion of these cases and more
references. For piecewise expanding maps of the interval, the result has
been obtained by A. Broise \cite{broise} under a generic non lacunarity
assumption. We show here that the result holds in full generality (as in
the independent case) for small enough values of $\alpha$. 
The proof is very similar to that of \cite{broise}. We give it below
for the sake of convenience.

\begin{lemma}\label{estime}
Under the assumptions of Theorem \ref{t3},
there is a compact neighborhood $K$ of the origin, and
two constants $0<c_{1}<c_{2}$ 
such that for any $\alpha\in K\backslash\{0\}$ and for any integer $k\geq
1+\beta^{-4}$ (recall that $\beta=\varphi'(\alpha)$)  we have
$$
\frac{c_{1}}{\beta\sqrt k}e^{-k\varphi(\alpha)}\le
\proba\big(S_{k}>k\alpha\big)
\le \frac{c_{2}}{\beta\sqrt k}e^{-k\varphi(\alpha)}\;.
$$

\end{lemma}

\proof

We first recall  a convenient representation of the probability we are
interested in.
For a given $\alpha\in K$, let $\beta=\varphi'(\alpha)$.
Define a sequence of positive measures $\nu_{k,\alpha}$ by
$$
d\nu_{k,\alpha}(s)=
\int_{0}^{1}dy\sum_{f^{k}(x)=y}\frac{e^{\beta
S_{k}(x)}e^{-kF(\beta)}}{\big|(f^{k})'(x)\big|}
\;h(x)\delta\big(s-k^{-1/2}(S_{k}(x)-k\alpha)\big)  \;,
$$
where $\delta$ denotes the Dirac measure. 
After a simple computation (see \cite{ney} or \cite{broise}) 
one gets
$$
\proba\big(S_{k}>k\alpha)=
e^{-k\varphi(\alpha)}\int_{0}^{\infty}e^{-\beta s\sqrt k}d\nu_{k,\alpha}(s)\;.
$$

In this expression of $\nu_{k,\alpha}$, one sees appearing  the  operator
$L_{z}$ given by
$$
\big(L_{z}v\big)(y)
=\sum_{f(x)=y}\frac{v(x)e^{zu(x)}}{\big|f'(x)\big|}\;.
$$

This operator depends analytically on $z$, and in the space of functions
of bounded variations, one can use analytic perturbation theory
\cite{kato}. In particular, there is a disk $D$ in the complex plane
centered at the origin such that for any $z$ inside that disk, the
operator $L_{z}$ has a peripheral spectrum consisting of a simple
eigenvalue denoted below $e^{F(z)}$. In the disk $D$, $F$ is
analytic and the corresponding eigenvector and eigenform depend also
analytically on $z$. Moreover the rest of the spectrum is inside a disk
of radius $\rho<1$ and there is a uniform bound on the corresponding
spectral projection. 

In particular, for a fixed real $\beta\in D$, 
 we have for any complex number $z$ such that $z+\beta\in D$
$$
\int_{-\infty}^{+\infty}e^{zs}d\nu_{k,\alpha}(s)= 
\int_{0}^{1}dy\sum_{f^{k}(x)=y}\frac{e^{\beta
S_{k}(x)}e^{-kF(\beta)}}{\big|{f^{k}}'(x)\big|}
\;h(x)e^{z(S_{k}(x)-k\alpha)/\sqrt k}
$$
$$
=e^{-kF(\beta)}e^{-\sqrt k\,\alpha z}
\int_{0}^{1}\big(L_{\beta+z/\sqrt k}\big)^{k}h(y)dy  \;.
$$
Note that in the above expression, since $u$ is bounded (and hence
$S_{k}$), the right-hand side is an entire function of $z$. We now apply the
analytic perturbation theory \cite{kato} in the last expression and get
$$
\int_{-\infty}^{+\infty}e^{zs}d\nu_{k,\alpha}(s)= 
C(\beta+z/\sqrt k)e^{-kF(\beta)+kF(\beta+z/\sqrt k)-k\sqrt\alpha}+
R_{k}(\beta+z/\sqrt k)
$$
where $R_{k}(w)$ is bounded for $w\in D$ by
$$
\big|R_{k}(w)\big|\le\Oun (3\rho/4)^{k}
$$
and $C(w)$ is analytic for $w\in D$.
Note that since $R_{k}$ can be written as the difference of two 
functions analytic in $D$, it is also analytic in $D$ and therefore we
have estimates of all its derivatives  in any smaller disk in terms of
its maximum in $D$. In particular, we have uniformly in $\alpha$ in a
neighborhood of the origin
$$
\int s\ d\nu_{k,\alpha}(s)= \Oun k^{-1/2}
$$
$$
\int s^{4}\ d\nu_{k,\alpha}(s)= \Oun 
$$
and
$$
\sigma(\alpha)^{2}
=\lim_{k\to\infty}\int s^{2}\ d\nu_{k,\alpha}(s)= F''(\beta)\;.
$$
As observed in \cite{broise}, if $F''(0)\neq 0$, that is to say $u$ is
not of the form $v-v\circ f$ for a function $v$ of bounded variation,
then by continuity we have $F''(\beta)\neq0$ in a neighborhood of the
origin in $\beta$ and also in a neighborhood of the origin in $\alpha$.

From these estimates, we can derive a Berry-Esseen estimate. 
We denote by $F_{k,\alpha}$ the law of $\nu_{k,\alpha}$ and by
$F_{\sigma}$ the law of the normal distribution with variance
$\sigma^{2}$. 

Under the assumptions of Theorem \ref{t1},
there is a number $\alpha'_{0}>0$ and a number $C>0$ such
that for any $|\alpha|<\alpha'_{0}$ we have
\begin{equation}\label{be}
\sup_{x}\big|F_{k,\alpha}(x)-F_{\sigma(\alpha)}(x)\big|\le \frac{C}{\sqrt k}\;.
\end{equation}

We refer to  \cite{feller} for the proof which uses a standard technique
once one has adequate control over the characteristic function. See also 
\cite{parrycoelho} and \cite{broise} for the case of dynamical systems. 

The main point of the above result is that the constant $C$ appearing
in (\ref{be}) is uniform in $\alpha$. 

To proceed with the proof, we will first give an asymptotic estimate of
$$
\int_{-\infty}^{+\infty}g(\beta s/\sqrt k)d\nu_{k,\alpha}(s)
$$
for fixed $C^1$ non-negative functions $g$ which tend to zero at
infinity and satisfying
$\int (s^2 g(s)+ |g'(s)|)\ ds<\infty$.
Integrating by parts we have
$$
\int_{-\infty}^{+\infty}g(\beta\sqrt k\,s)d\nu_{k,\alpha}(s)=
-\beta\sqrt k\int_{-\infty}^{+\infty}g'(\beta s\sqrt k)
F_{k,\alpha}(s) ds\;.
$$
(The same holds for the Gaussian measure in place of $\nu_{k,\alpha}$
and $F_\sigma$ in place of $F_{k,\alpha}$.)
Using (\ref{be}) we get uniformly for $\alpha$ in a neighborhood of
zero
$$
\left|\int_{-\infty}^{+\infty}g(\beta\sqrt k\,s)\ d\nu_{k,\alpha}(s)
-\int_{-\infty}^{+\infty}g(\beta\sqrt k\,s)\ d{\mathcal N}(0,\sigma^2 )(s)\right|=
$$
$$
\beta\sqrt k\ \left|\int_{-\infty}^{+\infty}\!\!\!\! g'(\beta s\sqrt k)
F_{k,\alpha}(s) ds
-\frac1{\sqrt{2\pi}\,\sigma(\alpha)}\int_{-\infty}^{+\infty}
\!\!\!\! ds\; g'(\beta s\sqrt k)
\int_{-\infty}^{s}\!\!\!\! e^{-\xi^{2}/(2\sigma(\alpha)^{2})}d\xi\right|
$$
$$
\le \Oun \beta \int_{-\infty}^{+\infty}
ds\; \big|g'(\beta s\sqrt k)\big|
\le\Oun k^{-1/2}\;.
$$
Using the hypothesis $k\geq 1+\beta^{-4}$, we obtain 
$$
\int_{-\infty}^{+\infty}g(\beta\sqrt k\,s)d\nu_{k,\alpha}(s)
=\frac{1}{\beta \sqrt{2\pi k}\,\sigma(\alpha)}\int_{-\infty}^{+\infty}
ds\; g(s)+\Oun k^{-1/2}\;.
$$
The main point of this estimate is that for $\beta$ small enough, the
first term on the right hand side dominates (see \cite{petrov} for the
exact coefficient  in the case of independent random variables).

The lemma follows by using $g(x)=e^{-\sqrt{1+x^2}}$ for an upper bound and 
for a lower bound by using a non negative $C^{1}$ function with compact
support in the interval $[1,2]$ and bounded above by $e^{-2}$.

\bigskip

The following lemma allows to control the probability of having
simultaneously $S_k>k\alpha$ and $S_k\circ f^r > k\alpha$.

\begin{lemma}\label{decouple}
Under the assumptions of Theorem \ref{t1},
there is a number $\alpha_{1}>0$ such that for any
compact subset $D$ of $(-\alpha_{1},0)\cup(0,\alpha_{1})$,
 there are two  constants $C>0$ and $c>0$
such that for any $\alpha\in D$, for any 
integers $k>0$ and $k>r>0$ we have
$$
\proba\big(S_{k}>k\alpha\;,\;S_{k}\circ f^{r}>k\alpha\big)\le
Ce^{-k\varphi(\alpha)-rc}\;.
$$
\end{lemma}

\proof We first observe that
$$
\proba\big(S_{k}>k\alpha\;,\;S_{k}\circ f^{r}>k\alpha\big)\le
\proba\big(S_{k}+\;S_{k}\circ f^{r}>2k\alpha\big)\;,
$$
and we will derive a large deviation estimate for this quantity.
Using the family of operators $L_{z}$ defined above, one gets
 for any real $\tilde\beta$
$$
\expectation\left( e^{\tilde \beta \big(S_{k}+\;S_{k}\circ f^{r}\big)}\right)
=\int  \big(L_{\tilde\beta}^{r}
L_{2\tilde\beta}^{k-r}L_{\tilde\beta}^{r}h\big)(x)\;dx\;.
$$
From the spectral theory of $L_{z}$, we get as in the proof of 
Lemma \ref{estime} uniformly for $\tilde\beta$ in a compact set
containing the origin in its interior
$$
\expectation \left(e^{\tilde \beta \big(S_{k}+\;S_{k}\circ f^{r}\big)}\right)
\le \Oun e^{2rF(\tilde\beta)}e^{(k-r)F(2\tilde\beta)}\;.
$$
We now use Chebychev inequality on the left hand side with
$\tilde\beta=\varphi'(\alpha)/2$ and obtain
$$
\proba\big(S_{k}+S_{k}\circ f^{r}>2k\alpha\big)\le
\Oun e^{-k\varphi(\alpha)-r(F(2\tilde\beta)-2F(\tilde\beta))}\;.
$$
From our assumptions on $u$ (in particular, $\sigma\neq0$), the function
$F$ is strictly convex, and therefore since $F(0)=0$ we have
$F(2\tilde\beta)-2F(\tilde\beta))>0$ for any $\tilde\beta\neq0$ in a
compact set. The lemma follows.

\section{Appendix: Strong approximation of return times by the measure of cylinders}\label{B}

\begin{lemma}\label{lemma}
Let $\epsilon>0$. Under the assumptions of Theorem \ref{t1} we have the following:
$$
-(1+\epsilon)\log n \leq \log\left[R_n(x)\mu(\mathcal{P}_n(x))\right]\leq \log\log(n^{1+\epsilon})
\quad\textup{eventually almost surely}\,.
$$
\end{lemma}

\proof
The lemma is the consequence of the exponential law for the distribution of the
random variables $R_n\mu(\mathcal{P}_n(\cdot))$. 
From \cite{paccaut-PhD} there exists a subset of cylinders
$\mathcal{P}_n^*\subset\mathcal{P}$ satisfying the following
conditions:
\begin{itemize}
\item 
$\zeta_n:=\mu\left(\mathcal{P}_n\backslash\mathcal{P}_n^*\right)$ 
with $\zeta_n\leq \Oun e^{-\kappa n}$, $\kappa>0$.

\item For any cylinder $C_n\in \mathcal{P}_n^*$, 
\begin{equation}\label{strong-approximation}
\sup_{t>0}\left|
\mu\left\{R_n\mu(C_n)>t\big|\ C_n\right\}-e^{-t}
\right|\leq c_1 e^{-c_2 n}
\end{equation}
where $c_1,c_2$ are positive constants,
$\mu\left\{\ \cdot\ \big|\ C_n\right\}$ denotes the conditional
expectation.

\end{itemize}

We want to find a summable upper-bound to 
$$
\mu\{\log\left[R_n \mu(C_n)\right]\geq \log t\}\leq
\sum_{C_n\in\mathcal{P}_n^*}\mu(C_n)
\mu\left\{\log\left[R_n \mu(C_n)\right]\geq \log t\big|\ C_n\right\}
+\zeta_n
$$
where $t$ will be chosen as a suitable sequence of positive real numbers.

Then, from (\ref{strong-approximation}), one gets for all $t>0$
$$
\mu\{\log[R_n \mu(C_n)]\geq \log t\}\leq \zeta_n+
\overbrace{\sum_{C_n\in\mathcal{P}_n^*}\mu(C_n)}^{\leq 1}\ (c_1 e^{-c_2 n}+ e^{-t})  \, .
$$
Take $t=t_n=\log(n^{1+\epsilon})$, $\epsilon>0$, to get
$$
\mu\{\log[R_n \mu(C_n)]\geq \log\log(n^{1+\epsilon})\}\leq  c_1 e^{-c_2 n} + \frac1{n^{1+\epsilon}}+\zeta_n\,.
$$
An application of the Borel-Cantelli lemma tells us that
$$
\log[R_n(x) \mu(\mathcal{P}_n(x))]\leq \log\log(n^{1+\epsilon})\quad\textup{eventually a.s.}\,.
$$
For the lower bound first observe that (\ref{strong-approximation}) gives, for all $t>0$
$$
\mu\{\log[R_n \mu(C_n)]\leq \log t\}\leq \zeta_n + c_1 e^{-c_2 n}+ 1-e^{-t}\leq \zeta_n + c_1 e^{-c_2 n}+ t\, .
$$
Choose $t=t_n=n^{-(1+\epsilon)}$, $\epsilon>0$, to get, proceeding as before,
$$
\log[R_n(x) \mu(\mathcal{P}_n(x))]\geq -(1+\epsilon)\log n\quad\textup{eventually a.s.}\,.
$$
This finishes the proof of the lemma.


\end{document}